\documentclass[11pt]{article}

\usepackage{color}
\usepackage[T1]{fontenc}
\usepackage{times}
\usepackage{graphics,graphicx}
\usepackage{amssymb}
\usepackage{subfigure}
\usepackage{epsfig}
\usepackage{amsbsy}
\usepackage{theorem}
\usepackage{float}
\usepackage{amsfonts}
\usepackage{natbib}
\usepackage{amssymb}
\usepackage{amsmath}
\usepackage{booktabs}

\newcommand{\beq}{\begin{equation}}
\newcommand{\eeq}{\end{equation}}

\newcommand{\beqar}{\begin{eqnarray}}
\newcommand{\eeqar}{\end{eqnarray}}

\begin{document}

\baselineskip=15pt

\title{{\bf A Framework Integrating the Dynamic Stiffness Matrix with Physics-Informed Neural Networks for Solving Eigenvalue Problems and Analysing Dynamic Response}}
\author{{\bf Yi An$^{(1)}$, Zhijiang Chen$^{(2,*)}$, Zhiqiang Feng$^{(3),(4)}$, Qian Cheng$^{(5)}$,} \\
{\bf Jack C.P. Cheng$^{(1)}$, Haijiang Li$^{(6)}$, Dalei Wang$^{(2,*)}$} \\ \\
$^{(1)}${\normalsize {\it Department of Civil and Environmental Engineering, }}\\
{\normalsize {\it The Hong Kong University of Science and Technology, Hong Kong, China.}} \\
$^{(2)}${\normalsize {\it College of Civil Engineering, Tongji University, Shanghai, 200092, China.}}\\
$^{(3)}${\normalsize {\it School of Mechanics and Aerospace Engineering,}}\\
{\normalsize {\it Southwest Jiaotong University, Chengdu, 610031, China.}} \\
$^{(4)}${\normalsize {\it LMEE Univ-Evry, Université Paris-Saclay, Evry, France.}}\\
$^{(5)}${\normalsize {\it College of Aerospace Engineering, Chongqing University, Chongqing, 400030, China.}}\\
$^{(6)}${\normalsize {\it School of Engineering, Cardiff University, The Parade, Cardiff, CF24 3AA, UK.}}\\
{\normalsize {\it  E-mail: 25512@tongji.edu.cn; wangdalei@tongji.edu.cn}} }
\maketitle

\begin{abstract}
This paper introduces a framework that integrates the dynamic stiffness matrix (DSM) with physics-informed neural networks (PINN). The DSM-PINN embeds physical constraints within the model and demonstrates robustness, particularly when addressing limited datasets across diverse investigations. In this approach, deep neural network outputs approximate the displacement fields of element nodes. Unlike the finite element method (FEM), the element shape functions are homogeneous solutions to the governing partial differential equation, forming the basis of the exact dynamic stiffness matrix, thereby avoiding high-order derivative terms. This matrix also serves as a frequency-domain spectral element, resulting in a strong-form PINN. The loss function is produced by connecting neural networks with dynamic stiffness matrices. We focus on utilising PINNs to resolve eigenvalue problems by employing the Wittrick-Williams algorithm, which overcomes the challenge of neural networks failing to converge to higher-order eigenvalues. Additionally, the frequency domain-PINN method is used to analyse structural dynamic responses under moving and impulsive loads, addressing the limitation of neural networks in handling complex numbers. Theoretical convergence stability of the suggested approach is also analysed even DSM is an indefinite matrix after implementing the boundary condition. The numerical results validate the practicality and efficacy of the recommended approach.

\end{abstract}

\section{Introduction}
The analysis of structural free vibration and dynamic response has always been a focal point of research, particularly in terms of stability and response characteristics. The research achievements in this field hold significant importance for industries such as civil engineering, aerospace, automotive, and shipbuilding, in terms of structural design optimization and safety performance enhancement. Over recent decades, numerous numerical computation techniques and derivatives of these techniques devised to address partial differential equations (PDEs) that characterise structural mechanical actions \citep{Samaniego1, LiangR1, BaiJ1}. Among them, the finite element method (FEM) \citep{LiHa1} is the most commonly utilised methodology, especially for problems in the low-to-medium frequency range. A primary benefit of the FEM is its capacity to accommodate intricate geometries. The Finite Element Method is fundamentally an approximation technique that depends on presumed displacement and coordinate functions (shape functions) and variational methods based on energy principles. This constrains its precision in the high-frequency spectrum, unless a substantial number of degrees of freedom and superior finite elements are utilised. Therefore, in the high-frequency range, the FEM becomes inefficient and higher cost for computation. Numerous engineering applications necessitate accurate free vibration analysis in the high-frequency spectrum, particularly for assessing the transmission of vibrational energy within structures (usually elastic wave energy). This aids in predicting issues such as fatigue damage and noise generation under high-frequency excitation, thereby ensuring the safety and comfort of structures.

Machine learning and data mining have become increasingly popular in a variety of fields, including materials science, computational science, and biology, as a direct result of the rapid development of information technology. In particular, neural network models, leveraging large-scale datasets for training and optimization, have achieved precise simulation of structural mechanical behaviors \citep{TsialiamanisG, LiangR1, HeW1}. This data-driven approach has demonstrated remarkable performance in areas involving structural health monitoring \citep{YeXW, DongY1, GaoY1, LiuJ1}, load identification \citep{YangH1}, uncertainty analysis \citep{LiuW1}, disaster management \citep{DengH1}, structural design \citep{XiangC2, VerduzcoL1}, and computational science \citep{BahmaniB1, MengX2, XuK2, Deshpande1, Krokos1}. Nonetheless, while these breakthroughs and considerable opportunity, contemporary machine learning techniques often encounter difficulties in deriving interpretable insights from data. Moreover, models that only use data frequently struggle with extrapolating and actual biases, which can result in estimations that are neither reliable or physically compatible. These limitations eventually hinder the generalization potential of artificial neural networks \citep{Karniadakis1}. Essentially, machine learning aims to discover hidden manifold structures within data samples and the probability distributions defined on them. However, the question of how to leverage mechanical principles and engineering logic to determine parameterised representations of coding mappings with stronger learning capabilities and better stability remains an urgent area for in-depth exploration.

Recent research improvement has focused on inclusion of physical principles and domain information into deep learning models, with the objective of addressing the shortcomings of exclusively data-driven models \citep{RaissiM1, Weinan1, Samaniego1, YinM1, ChiuP1, MengX1}. This strategy, referred to as physics-informed deep learning, augments the ability to forecast of machine learning models by incorporating observational information, physical rules, and empirical findings. It utilises physical principles to constrain data and empirical knowledge to provide prior information, thereby improving the model's predictive ability \citep{CuomoS1, BaiJ1}. A notable progression in this domain is the creation of Physics-Informed Neural Networks (PINNs) \citep{RaissiM1}. PINNs include governing partial differential equations (PDEs) inside deep neural networks, achieving the integration of physical knowledge and data. Based on this, a strong-form loss function is formulated and used as an inductive bias (inductive bias refers to the assumptions made by machine learning algorithms during training, which affect the algorithm's data fitting method and generalization ability) \citep{LuL1}. This method improves the interpretability of neural networks by imposing physical constraints. Furthermore, even in the presence of incomplete or unlabeled data, PINNs exhibit robust predictive performance because they consistently generate outputs that conform to physical laws, thereby enhancing generalization and extrapolation capabilities.

Physical Information Neural Networks tend to be classified into two primary categories according to their theoretical foundation and solution approach, each with different variants: governing equation-based PINNs (strong-form) \citep{RaissiM1, Amini1, Jagtap1} and energy-based PINNs (weak-form) \citep{Samaniego1, Weinan1, LiHa1, XiongW1, HuangM1}. Governing equation-based PINNs typically operate at collocation points, using automatic differentiation to obtain the residuals of PDEs to drive the minimization of these residuals \citep{HeW1, LiangR1, MaoZ1}. Notable methods in this category encompass the Deep Galerkin Method (DGM) \citep{Sirignano1}, the hp-Variational PINN (hp-VPINN) \citep{Kharazmi1}, the Conservative PINN (CPINN) \citep{Jagtap1}, and the Physics-Informed Graph Neural Networks (PIGCN) \citep{GaoH1}. Conversely, energy-based PINNs solve issues by finding the extreme point of an energy value as a loss function, reducing this functional yields the solution to the PDE \citep{Samaniego1}. The Deep Ritz Method (DRM) \citep{Weinan1} is among the pioneering techniques to establish this framework of principles. Thereafter, the Deep Energy Method (DEM) \citep{Nguyen1} was utilised in computational mechanics, and its effectiveness in solving elasticity problems involving cracks and inclusions was demonstrated by \citep{WangY4}. Furthermore, inspired by the Finite Element Method and PINNs, within the framework of Energy-based PINNs, the linear equation system in FEM can be incorporated into the loss function, providing stability, convergence, and practicality for complex engineering problems \citep{XiongW1}. Energy-based PINNs have advantages such as simplified loss functions and diminished derivative orders, rendering them typically more facile to train than governing differential equation-based PINNs.
 
Due to the insufficient similarity between the neural network output and the corresponding high-order feature functions during training, it may fail to converge to high-order feature modes and often converges to lower-order feature modes \citep{YooS}. Furthermore, existing research indicate that neural networks encounter difficulties in modeling functions with multi-frequency attributes \citep{LiangR1, TancikM1}, which are prevalent in structural dynamics. This limitation, referred to spectral bias \citep{Fridovich-KeilS}, restricts the application of neural networks in such environments. To tackle these challenges, this research presents an innovative framework combining the Dynamic Stiffness Matrix (DSM) with the Physics-Informed Neural Network (PINN). Since extensively recorded in the literature concerning the benefits of DSM \citep{Wittrick1, Banerjee3, Yuan2}, DSM has significant advantages over the Finite Element Method (FEM) that relies on approximate interpolation functions and energy variational principles. It utilises frequency-dependent components derived from the homogeneous solutions of free vibration governing equations. The DSM has proven widely employed in multiple fields through the Wittrick-Williams algorithm, including beam analysis \citep{Liux2, Banerjee9, ChenZ2, ChenZ3}; multi-cracked beam or frame studies \citep{CaddemiS1, HanF, CannizzaroF1}; plate analysis \citep{Banerjee2, Papkov1, Liux3, Chauhan1}; shell analysis \citep{Nihal1}; composite material analysis \citep{Shams1}; CUF element analysis \citep{Liu4}; lattice structures \citep{Liux1}; and membrane structures \citep{Kim1}. The high accuracy of DSM across the entire frequency range makes it particularly valuable in structural buckling and vibration analysis, both of which involve eigenvalue problems. Furthermore, DSM markedly diminishes the total degrees of freedom necessary for computation and cuts the corresponding computational expenses. The spectral element method in the frequency domain is also considered an accurate element method that provides precise solutions for structural dynamics problems in the frequency domain \citep{Liu5, KimT1, KimT2, KimT3, KimT4}. In fact, spectral elements are equivalent to dynamic stiffness elements \citep{LeeU}. Therefore, a frequency-domain PINN technique is developed inside the DSM framework to mitigate the spectral bias problem of neural networks.

The primary advances of the suggested method are as follows. First, A deep learning system is introduced that incorporates the dynamic stiffness method (DSM) into neural networks. The eigenmode of an arbitrary-order structure is reformulated as a neural network lowest-eigenvalue problem using the Wittrick-Williams algorithm. Second, a robust and integrated framework for frequency-domain physics-informed neural networks is developed, which combines experimental data with physical laws. The method applies the Fourier transform to translate the original time domain dynamic load into the frequency domain, thereby simplifying neural network simulations and enhancing accuracy. Additionally, a stable loss function is developed following the properties of the dynamic stiffness matrix, facilitating neural network training and enabling a unified strategy for both simple and complicated applications. A series of case studies examining the dynamic response of beams under diverse boundary conditions is performed to thoroughly assess the method's efficacy. The following section of the paper is structured as outlined below. Section 2 briefly introduced the fundamentals of stiffness-matrix-type PINNs and their calculation. Section 3 delineates the methodology for transforming the eigenmode of arbitrary order into a lowest-eigenvalue problem and provides numerical results for beam and frame problems. Section 4 describes the proposed frequency-domain PINN method. To ascertain the efficacy of the PINN, the implicit time-integration method and the Newmark-$\beta$ method will be used. Ultimately, Section 5 ends the study and offers a perspective for future endeavors.

\section{The physics-informed neural networks and their application in finite element method}

Since the complete dynamic stiffness matrix is assembled using a method analogous to the finite element method, this section presents a comprehensive method for embedding the stiffness matrix into the loss function of artificial neural networks. This allows for the integration of physics-informed neural networks with FEM. In this section, the overall methodology, boundary conditions (BCs), and algorithmic implementation are all covered. This strategy can be utilised for other partial differential equations (PDEs), which are characterised by the presence of a system of linear equations derived from energy variational formulations.

\subsection{Energy- and Equilibrium Equation-Based PINNs in FEM}
Boundary value problems discretised via FEM can be solved using energy-based PINNs \citep{XiongW1, BastekJ}, which are known as deep finite element methods. The method partitions the solution domain into discrete elements, guaranteeing that the test functions within each one retain a uniform and simple form. Furthermore, it approximates the solution domain of nodal variables by employing neural networks $\mathcal{N}(\textbf{x}; \theta)$:
 \begin{equation}
  \label{eq4}
\textbf{u}(\textbf{x}_j)=\mathcal{N}(\textbf{x}_j; \theta), \ \ j=1,2,...n,
 \end{equation}
where $\textbf{u}(\textbf{x}_j)$ denotes the displacement vector corresponding to the $j$-th node, where $n$ signifies the total count of nodes, and $\textbf{x}_j$ indicates the global coordinates associated with the $j$-th node. The neural network parameters $\theta$ are composed of weight matrices and bias vectors. Then, the energy-based FEM can be derived from the potential energy:
  \begin{equation}
  \label{eq7}
\Pi=\frac{1}{2}{\textbf{u}(\textbf{x})}^{T}\cdot{\textbf{K}}\cdot\textbf{u}(\textbf{x})-{\textbf{u}(\textbf{x})}^{T}\cdot\textbf{P},
 \end{equation}
where $\textbf{K}$ denotes the global stiffness matrix, and $\textbf{P}$ indicates the global nodal load vector. The global nodal displacement vector $\textbf{u}(\textbf{x})$ is is capable of being approximated using artificial neural networks $\mathcal{N}(\textbf{x}; \theta)$. Therefore, the objective is to ascertain the neural network parameters $\theta^*$ that optimise the potential energy:
  \begin{equation}
  \label{eq11}
\theta^*=arg\min_\theta \left(\mathcal{L}_{phys}=\frac{1}{2}{\mathcal{N}(\textbf{x}; \theta)}^{T}\cdot{\textbf{K}}\cdot\mathcal{N}(\textbf{x}; \theta)-{\mathcal{N}(\textbf{x}; \theta)}^{T}\cdot\textbf{P}\right).
 \end{equation}
 
There is another method to embed the physical constraints in PINNs by employing the finite element method. This method is comparable to the governing equations-based PINNs (strong form) \citep{RaissiM1}. The optimisation of the PINNs' loss function $\mathcal{L}_{phys}$ is denoted \citep{HuangM1, BastekJ}:
  \begin{equation}
  \label{eq13}
\theta^*=arg\min_\theta \left(\mathcal{L}_{phys}=\frac{1}{n}\sum^n_{i=1}\left[\textbf{K}_{ij}\cdot{\mathcal{N}(\textbf{x}_j; \theta)}-\textbf{P}_i\right]^2\right),
 \end{equation}
where $\textbf{K}$, $\textbf{P}$, and $\mathcal{N}(\textbf{x}; \theta)$ remain the global stiffness matrix, the global nodal load vector, and the outputs of the deep learning model. In both energy-based and equilibrium-based PINNs, the data loss, which quantifies the model's fit to observed or experimental data, is defined as follows:
  \begin{equation}
  \label{eq14}
\mathcal{L}_{dat}=\alpha \sum_k\|\mathcal{N}(\textbf{x}_k; \theta)-\hat{\textbf{u}}(\textbf{x}_k)\|^2,
 \end{equation}
where $\hat{\textbf{u}}(\textbf{x}_k)$ represents the data from experiments at the location $\textbf{x}_k$, and $\alpha$ serves as the penalty parameter. The complete loss function, encompassing both physical and measurable data components, can be stated as follows:
  \begin{equation}
  \label{eq15}
\mathcal{L}=\mathcal{L}_{phys}+\mathcal{L}_{dat}.
 \end{equation}

Resolving this system is equivalent to determining the optimal set of neural network parameters $\theta^*$ that minimise the entire loss function:
  \begin{equation}
  \label{eq16}
\theta^*=arg\min_\theta \left(\mathcal{L}_{phys}+\mathcal{L}_{dat}\right).
 \end{equation}
 
Optimisation strategies in machine learning primarily consist of variants of first-order methods, including the gradient descent method \citep{Sirignano1}. If only physical loss were considered, the optimisation algorithm would iterate as follows:
  \begin{equation}
  \label{eq12}
\theta_{n+1}=\theta_{n}-\gamma_n\nabla_\theta(\mathcal{L}_{phys}),
 \end{equation}
where $\theta_{n}$ denotes the parameters of the neural network being trained at the present iteration, and $\theta_{n+1}$ specifies the parameters at the subsequent iteration. The parameter $\gamma_n$ is the learning rate, and $\nabla_\theta(\mathcal{L}_{phys})$ represents the first-order gradient. Various optimisation algorithms can be employed, such as stochastic gradient descent (SGD), Adaptive Moment Estimation (Adam), and the quasi-Newton method L-BFGS.

\subsection{Boundary conditions}
\begin{figure}[t]
  \centering
  \includegraphics[scale=0.88]{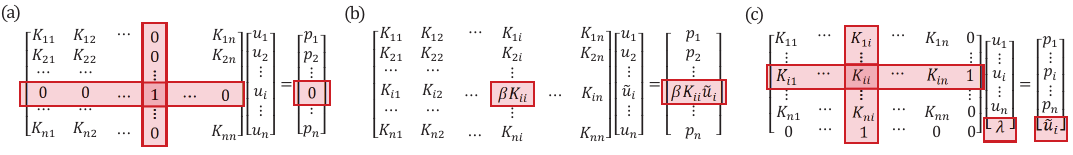}
\caption{The three ways to enforce boundary conditions in the FEM.}
\label{fig1}
\end{figure}

Boundary conditions (BCs) are typically categorised as either natural or essential. Natural BCs are integrated via the nodal load vector $\textbf{P}$. Essential BCs can be introduced as observable data, as shown in Eq. \eqref{eq14}, within the loss function term, and are referred to as soft BCs \citep{RaissiM1, Sirignano1}. Alternatively, hard BCs can be implemented in PINNs by customising the neural network architecture to inherently satisfy the BCs \citep{BergJ, XiongW1}.
  \begin{equation}
  \label{eq17}
\mathcal{N}(\textbf{x}; \theta)=\textbf{D}(\textbf{x})\tilde{\mathcal{N}}(\textbf{x}; \theta)+\textbf{G}(\textbf{x}),
 \end{equation}
where $\textbf{D}(\textbf{x})$ signifies the distance function, which is zero on the essential BCs, whereas $\textbf{G}(\textbf{x})$ represents the seamless extension of the essential BCs. Both $\textbf{D}(\textbf{x})$ and $\textbf{G}(\textbf{x})$ could be defined as either manually crafted functions or neural networks. Based on prior research \citep{HuangM1, BastekJ, XiongW1}, BCs can be directly applied to the outputs of PINNs, similar to conventional finite element methods, the framework for implementing BCs on the stiffness matrix and loss function as depicted in Fig. \ref{fig1} (a). The global stiffness matrix $\textbf{K}$ is altered by assigning a value of one to the main diagonal elements and a value of zero to all the other components in the corresponding rows and columns. The relevant components in the nodal force vector are concurrently assigned a value of $0$. In cases where a specified nodal displacement is a nonzero value $\tilde{u}_i$, the diagonal element $\textbf{K}_{ii}$ in the stiffness matrix $\textbf{K}$ is multiplied by a large penalty parameter $\beta$, and $\beta\textbf{K}_{ii}\tilde{u}_i$ replaces $\textbf{P}_i$ in the nodal load vector, a technique known as the penalty method (see Fig. \ref{fig1} (b)). Alternatively, the Lagrange multiplier method may be employed, as depicted in Fig. \ref{fig1} (c).
The following equation does not hold because the BCs are not incorporated into the system of linear equations (usually, a stiffness matrix without sufficient BCs is a singular matrix, and the nodal force vector does not include the reaction force components),
\begin{equation}
  \label{eqequal}
\textbf{K}\cdot\mathcal{N}(\textbf{x}; \theta)=\textbf{P},
 \end{equation}
but is satisfied by the modified global stiffness matrix $\textbf{K}^*$ and nodal force vector $\textbf{P}^*$. The matrix $\textbf{K}^*$ is presently symmetric and positive definite.
\begin{equation}
  \label{eq18}
\textbf{K}^*\cdot\mathcal{N}(\textbf{x}; \theta)=\textbf{P}^*.
 \end{equation}
Eq. \eqref{eq18}, when combined with either loss function \eqref{eq15}, directly enforces BCs on the outputs of the deep learning model. This approach remains compatible with both soft and hard boundary conditions.

\subsection{The second-order derivative of loss function}
To accelerate the convergence of the loss function, neural network optimisation methods have been developed. At each step of the learning process, the first-order gradient of the loss function with respect to the parameters $\theta$ of the neural network progressively approaches zero. On the one hand, the first-order gradient of the energy-based PINNs is delineated in \citep{XiongW1} within the FEM loss function \eqref{eq11}:
\begin{equation}
  \label{eq19}
\frac{\partial\mathcal{L}}{\partial \theta}=(\textbf{K}^*\cdot\mathcal{N}(\textbf{x}; \theta)-\textbf{P}^*)^T\cdot\frac{\partial\mathcal{N}}{\partial \theta}.
 \end{equation}
The expression $\textbf{K}^*\cdot\mathcal{N}(\textbf{x}; \theta)-\textbf{P}^*$ pertains to solutions obtained from the FEM. This relationship illustrates that, in energy-based PINNs utilised in finite element methods, the solution completely fulfills the partial differential equations. The gradient, as indicated by Eq. \eqref{eq19}, tends toward zero as the neural network model reaches an extremum. The second-order derivative (Hessian matrix) of the loss function \eqref{eq11} with respect to the neural network parameters is expressed as:
\begin{equation}
  \label{eq20}
\mathbb{H}_{\mathcal{L}}=\frac{\partial^2\mathcal{L}}{\partial \theta^2}=\frac{\partial\mathcal{N}}{\partial \theta}^T\cdot\textbf{K}^*\cdot\frac{\partial\mathcal{N}}{\partial \theta}+(\textbf{K}^*\cdot\mathcal{N}(\textbf{x}; \theta)-\textbf{P}^*)^T\cdot\frac{\partial^2\mathcal{N}}{\partial \theta^2}.
 \end{equation}
In the above equation, $(\textbf{K}^*\cdot\mathcal{N}(\textbf{x}; \theta)-\textbf{P}^*)^T\frac{\partial^2\mathcal{N}}{\partial \theta^2}$ is the zero matrix when the model converges. Moreover, the stiffness matrix $\textbf{K}$ in FEM is a positive semi-definite matrix and $\textbf{K}^*$ is a positive matrix with sufficient constraints. This indicates that energy-based PINNs in FEM possess a degree of local stability. 

On the other hand, the loss function for the nodal equilibrium equations based on PINNs in FEM \eqref{eq13} can be transformed into the following expressions:
\begin{equation}
  \label{eq21}
\mathcal{L}=\frac{1}{n}(\textbf{K}^*\cdot\mathcal{N}(\textbf{x}; \theta)-\textbf{P}^*)^T\cdot(\textbf{K}^*\cdot\mathcal{N}(\textbf{x}; \theta)-\textbf{P}^*)=\frac{1}{n}(\mathcal{N}^T\cdot(\textbf{K}^*)^2\cdot\mathcal{N}-2\mathcal{N}^T\cdot(\textbf{K}^*)^T\cdot P^*+(P^*)^T\cdot P^*).
 \end{equation}
The stiffness matrix (global or elemental) in the finite element method is symmetric, such that $\textbf{K}^T=\textbf{K}$. The first-order gradient of the nodal equilibrium equations, formulated using a deep learning model in the FEM-type loss function \eqref{eq21}, with respect to the parameters $\theta$ of the neural network, is given by:
 \begin{equation}
  \label{eq22}
\frac{\partial\mathcal{L}}{\partial \theta}=\frac{2}{n}((\textbf{K}^*)^2\cdot\mathcal{N}-\textbf{K}^*\cdot P^*)^T\cdot\frac{\partial\mathcal{N}}{\partial \theta}.
 \end{equation}
The term $(\textbf{K}^*)^2\cdot\mathcal{N}-\textbf{K}^*\cdot P^*$ continues to represent the solutions of the finite element method. Consequently, the gradient described by Eq. \eqref{eq21} must be zero as the neural network model converges to the extremum. The Hessian matrix of the loss function in Eq. \eqref{eq21} concerning the parameters $\theta$ of the neural network, is given as follows:
 \begin{equation}
  \label{eq23}
\mathbb{H}_{\mathcal{L}}=\frac{\partial^2\mathcal{L}}{\partial \theta^2}=\frac{2}{n}\frac{\partial\mathcal{N}}{\partial \theta}^T\cdot(\textbf{K}^*)^2\cdot\frac{\partial\mathcal{N}}{\partial \theta}+\frac{2}{n}((\textbf{K}^*)^2\cdot\mathcal{N}(\textbf{x}; \theta)-\textbf{K}^*\cdot\textbf{P}^*)^T\cdot\frac{\partial^2\mathcal{N}}{\partial \theta^2}.
 \end{equation}
Similar to energy-based PINNs in FEM, the $((\textbf{K}^*)^2\cdot\mathcal{N}-\textbf{K}^*\cdot P^*)^T\cdot\frac{\partial^2\mathcal{N}}{\partial \theta^2}$ remain zero matrix as the model converges. $(\textbf{K}^*)^2$ remain a positive semi-definite matrix or a positive matrix with sufficient constraints. This demonstrates the nodal equilibrium equations based on PINNs in FEM, which exhibit a degree of local stability. 

Optimization of the loss function Eq. \eqref{eq15}, incorporates both physics-based rules (FEM) and data-driven losses to determine the optimal parameters $\theta$. This approach ensures that the model satisfies the governing equations and aligns with observed data. Regardless of whether the nodal equilibrium equations are formulated using standard PINNs or energy-based PINNs, the stiffness matrix is consistently derived from polynomial shape functions to approximate the solution of the governing partial differential equation through energy variation (principle of minimum potential, Hu-Washizu principle, Galerkin method). Therefore, both approaches are classified as weak-form methods.

\section{Addressing the eigenvalue problems in the dynamic stiffness matrix using PINNs}
\subsection{Dynamic stiffness matrix}
This section introduces the dynamic stiffness matrix (DSM) using the Euler-Bernoulli beam model. In addition, the procedures for DSM-PINNs can be utilised for other structures as introduced in the earlier work because the key steps are only associated with DSM. The length of the beam is denoted by $L$, its width is denoted by $b$, its height is denoted by $h$, and its mass per unit length is denoted by $m$. The transverse deflection of the beam is indicated by $w(x,t)$, where $x$ represents the coordinate along the beam length. The equation that governs free vibration in the absence of damping is a differential equation of the fourth order:
 \begin{equation}
  \label{eq24}
EI\frac{\partial^4 w(x,t)}{\partial x^4}+m\frac{\partial^2 w(x,t)}{\partial t^2}=0.
 \end{equation}
The solution to Eq. \eqref{eq24} should be obtained by employing the method of separation of variables:
 \begin{equation}
  \label{eq25}
w(x,t)=\bar{w}(x)e^{\textmd{i}\omega t}.
 \end{equation}
The amplitude of transverse deflection $w$ and the circular test frequency are represented by $\bar{w}$ and $\omega$, respectively. Consequently, Eq. \eqref{eq24} becomes the governing equation for the transverse deflection amplitude $\bar{w}$:
 \begin{equation}
  \label{eq26}
EI\frac{d^4 \bar{w}}{d x^4}-m\omega^2\bar{w}=0.
 \end{equation}
In fact, the Eq. \eqref{eq26} can be seen as the Fourier transform of Eq. \eqref{eq24} rather than assuming $w(x,t)=\bar{w}(x)e^{\textmd{i}\omega t}$ but $\bar{w}=\int w(x,t)e^{\textmd{i}\omega t}dt$. The governing equation becomes the same form, similar to Eq. \eqref{eq25}, because the discrete Fourier transform is used: 
 \begin{equation}
  \label{eqdiscrete}
w(x,t)=\frac{1}{N}\sum^{N-1}_{k=0}\bar{w}(x,\omega_k)e^{\frac{2ikn\pi}{N}}
  \end{equation}
The homogenous solutions for the ordinary differential equation \eqref{eq26}:
 \begin{equation}
  \label{eq27}
\bar{w}(\zeta, \omega)=A_1\sin \alpha \zeta +A_2 \cos \alpha \zeta + A_3 \sinh \alpha \zeta +A_4 \cosh \alpha \zeta 
 \end{equation}
where $\alpha=(\frac{m\omega^2L^4}{EI})^{1/4}$ and $\zeta=x/L$. Therefore, a connection between the amplitudes of the edge bending moment $\bar{M}$, shear force $\bar{F}_s$, deflection $\bar{w}$, and rotation $\bar{\phi}$ can be expressed as:
 \begin{equation}
  \label{eq28}
\bar{\phi}=\frac{d\bar{w}}{dx}, \ \ \bar{M}=EI\frac{d^2\bar{w}}{dx^2}, \ \ \bar{F}_s=EI\frac{d^3\bar{w}}{dx^3}.
 \end{equation}
Therefore, the system of linear equations for displacement can be rewritten as:
 \begin{equation}
   \label{eq29}
\Delta=\textbf{R}\textbf{A},
 \end{equation}
where $\Delta=[\bar{w}(0),\bar{\phi}(0),\bar{w}(1),\bar{\phi}(1)]^T$ and $\textbf{A}=[A_1, A_2, A_3, A_4]^T$. In a similar manner, the relationship between $\textbf{P}=[-\bar{F}_s(0),\bar{M}(0),\bar{F}_s(1),-\bar{M}(1)]^T$ and the constant vector $\textbf{A}$, along with the corresponding type of homogeneous solutions, is expressed as follows:
 \begin{equation}
   \label{eq30}
\textbf{P}=\textbf{U}\textbf{A}.
 \end{equation}
The constant vector $\textbf{A}$ may be eliminated from Eqs \eqref{eq29} and \eqref{eq30} to derive the $4\times4$ DSM based the Euler-Bernoulli beam theory. This process involves inverting the square matrix $\textbf{R}$, followed by multiplication by the square matrix $\textbf{U}$.
 \begin{equation}
   \label{eq31}
\textbf{P}=\textbf{K}\Delta, \ \ \textbf{W}=\textbf{U}\textbf{R}^{-1}.
 \end{equation}
Because the shape functions are derived from a homogeneous solution, the matrix $\textbf{R}$ and $\textbf{U}$ are `exact'. Consequently, the numerical computation produces the dynamic stiffness matrix $\textbf{W}$ without employing energy variation and the results should converge to the analytical solution. The matrix $\textbf{R}$, $\textbf{U}$ and $\textbf{W}$ are shown in the Appendix. Thus, this dynamic stiffness matrix, combined with nodal equilibrium equations derived from PINNs, yields a strong-form solution for structural problems. The global DSM is constructed utilising a procedure similar to the FEM, establishing a foundation for accurate structural eigenvalue calculations with related algorithms. The eigenvalue problem is solved through the Wittrick-Williams algorithm \citep{Wittrick1} and Newton's method \citep{Yuan2}, both recognised for their accuracy and reliability. By selecting continuous test frequencies, the method accurately brackets any analytical natural frequency within selected upper and lower bounds. As the Wittrick-Williams technique confines natural frequencies to a narrow range, PINNs are able to converge and generate the corresponding mode shapes.

\subsection{Solving the eigenvalue problem with DSM-PINNs}
Inspired by the idea of the second-order mode-finding method in DSM \citep{Yuan2}, the eigenvalue problem becomes a homogeneous linear equation after using the Wittrick-Williams algorithm:
\begin{equation}
 \label{eq34} 
\textbf{W}^*(\omega)\Delta=0, \ \ \omega\in(\omega_l, \omega_u),
 \end{equation} 
where $\omega_l$ and $\omega_u$ are the lower and upper circular frequency bounds. $\textbf{W}^*$ is the DSM after imposing the BCs as described in Section 2. Assuming the $j$-th order natural frequency is to be determined, the following expressions can be derived:
\begin{equation}
 \label{eq35} 
j-1=j_0(\omega_l)+s\{\textbf{W}^*(\omega_l)\}, \ \ j=j_0(\omega_u)+s\{\textbf{W}^*(\omega_u)\}
 \end{equation} 
where $j_0$ denotes the entire count of natural frequencies of the structure with all degrees of freedom clamped. Although the DSM method is an accurate algorithm, it only converges the natural frequency to a very small interval $(\omega_l, \omega_u)$. Therefore, with $\omega_a$ denoting an approximation to the exact $j$-th order eigenvalue $\omega_e$, a second-order Taylor expansion is performed on the DSM:
\begin{equation}
 \label{eq36} 
\textbf{W}^*(\omega_e)=\textbf{W}^*(\omega_a)+(\omega_e-\omega_a)\frac{d\textbf{W}^*(\omega_a)}{d\omega}+\frac{(\omega_e-\omega_a)^2}{2}\frac{d^2\textbf{W}^*(\omega_a)}{d\omega^2}+o((\omega_e-\omega_a)^3).
 \end{equation} 
The exact eigenvector $\textbf{d}_e$ can be substituted into the Taylor expansion \eqref{eq36}, ignoring the second derivatives:
\begin{equation}
 \label{eq37} 
\textbf{W}^*(\omega_e)\textbf{d}_e\approx\textbf{W}^*(\omega_a)\textbf{d}_e+(\omega_e-\omega_a)\frac{d\textbf{W}^*(\omega_a)}{d\omega}\textbf{d}_e.
 \end{equation} 
Since $\textbf{W}^*(\omega_e)\textbf{d}_e=0$ by Eq. \eqref{eq34}, the expression \eqref{eq37} becomes:
\begin{equation}
 \label{eq38} 
\textbf{W}^*(\omega_a)\textbf{d}_e\approx-(\omega_e-\omega_a)\frac{d\textbf{W}^*(\omega_a)}{d\omega}\textbf{d}_e.
 \end{equation} 
The $\omega_a$ can be the centre of the frequency range between $\omega_l$ and $\omega_u$, so $\omega_a=(\omega_l+\omega_u)/2$ and the derivative of the DSM is approximated by the central difference method below:
\begin{equation}
 \label{eq39} 
\frac{d\textbf{W}^*(\omega_a)}{d\omega}\approx\frac{\textbf{W}^*(\omega_u)-\textbf{W}^*(\omega_l)}{\omega_u-\omega_l}.
 \end{equation}
Eq. \eqref{eq39} is an excellent approximation because it only has one eigenvalue within the small frequency range according to the Wittrick-Williams algorithm. Then, substituting Eq. \eqref{eq39} into \eqref{eq38}, a standard generalised eigenproblem can be obtained:
\begin{equation}
 \label{eq40} 
\left(\textbf{W}^*(\omega_a)-\bar{\mu}(\textbf{W}^*(\omega_u)-\textbf{W}^*(\omega_l))\right)\textbf{d}_e=0,
 \end{equation}
where the generalised eigenvalue $\bar{\mu}$ is the lowest eigenvalue of the system \eqref{eq40} \citep{Yuan2}. There are two key challengs when using PINNs to solve eigenvalues: Firstly, the governing equation (free vibration or buckling) is a homogeneous differential equation, so the non-trivial solutions of the eigenmode solutions need to be solved \citep{HarcombeL, YangQ1} and normalised. Secondly, the system tends to converge to low-order eigenmodes during PINNs training, which is related to spectral bias \citep{YooS}. Spectral bias describes the tendency of neural networks to converge to locally optimal solutions with low frequencies. However, with the help of the Wittrick-Williams algorithm, the arbitrary-order eigenvalue has been transformed into the lowest one $\bar{\mu}$. Then, the $\bar{\mu}$ can be expressed as:
\begin{equation}
 \label{eq41} 
\bar{\mu}\approx\frac{\omega_a-\omega_{\bar{\mu}}}{\omega_u-\omega_l}.
 \end{equation}
The eigenvalue $\bar{\mu}$ can be checked by Eq. \eqref{eq41} with the condition $|\bar{\mu}|<1/2$. Hence, if the eigenvalues $\bar{\mu}$ is determined, a better eigenvalue $\omega_{\bar{\mu}}$ can be obtained:
\begin{equation}
 \label{eq42} 
\omega_{\bar{\mu}}=\omega_a-\bar{\mu}(\omega_u-\omega_l).
 \end{equation}

\begin{figure}[htb]
  \centering
  \includegraphics[scale=0.8]{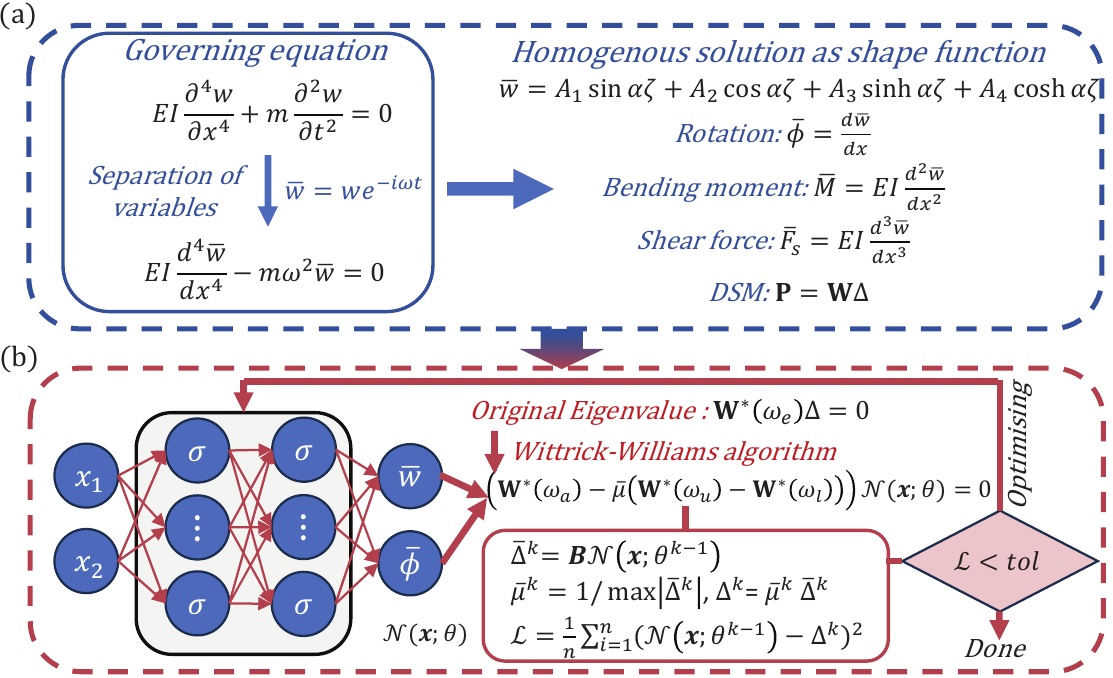}
\caption{Solution procedure for the structural eigenvalue problem using the DSM-PINNs method. (a). The dynamic stiffness matrix (DSM) formulation is represented, with the Euler-Bernoulli beam as an example. (b) The eigenvalue for arbitrary order can be transformed into the lowest generalised eigenvalue problem by the Wittrick-Williams algorithm. The neural network is used to solve the eigenvector. Note that this method can be extended to other structural elements.}
\label{fig3}
\end{figure}

Inspired by inverse iteration, the neural network $\mathcal{N}$ approximates the eigenvector $\Delta^{k}$ at the $k$-th iteration during training. The inverse iteration procedure, which is guaranteed to converge to the eigenpair (eigenvalue and associated eigenvector) corresponding to the target eigenvalue, is analogous to the inverse power method neural network \citep{YangQ1}:
\begin{equation*}
\bar{\Delta}^{k}=\textbf{B}\cdot\mathcal{N}(\textbf{x}; \theta^{k-1}),
 \end{equation*}
\begin{equation*}
\bar{\mu}^{k}=\frac{1}{\delta^{k}}, \ \textmd{with} \ \delta^{k}=\max|\bar{\Delta}^{k}| \ \textmd{and} \ {\Delta}^{k}=\bar{\mu}^{k}\bar{\Delta}^{k},
 \end{equation*}
\begin{equation}
 \label{eq43} 
\mathcal{L}=\frac{1}{n}\sum^n_{i=1}(\mathcal{N}(\textbf{x}; \theta^{k-1})-{\Delta}^{k})^2,
 \end{equation}
where matrix $\textbf{B}$ is that:
\begin{equation}
 \label{eq44} 
\textbf{B}={(\textbf{W}^*(\omega_a))}^{-1}(\textbf{W}^*(\omega_u)-\textbf{W}^*(\omega_l)).
 \end{equation}
Upon convergence of the neural network, the smallest eigenvalue $\bar{\mu}$ and the corresponding eigenvector, represented by $\mathcal{N}(\textbf{x}; \theta)$, are determined. Because the eigenvector is normalised in each training iteration (as expressed in Eq. \eqref{eq43}), the nontrivial solution can be obtained. In addition, Eq. \eqref{eq41} shows that the eigenvalue problem for an arbitrary $j$-th order mode is transformed into the lowest order eigenvalue problem. Thus, the two difficulties in using neural networks to solve the eigenvalue problem have been addressed. The network structure is shown in Fig. \ref{fig3}. The following subsection presents case studies on solving eigenproblems for beam and frame structures (natural frequencies) using DSM-PINNs.

\subsection{Case studies for eigenvalue problems}
\begin{figure}[htb]
  \centering
  \includegraphics[scale=0.6]{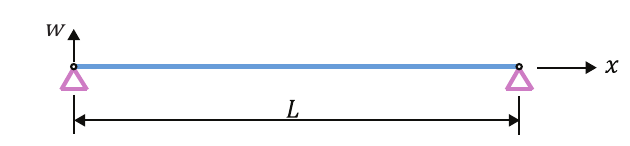}
\caption{The sketch of beam structures with simply supported boundary conditions.}
\label{fig4}
\end{figure}

To ensure robust validation, results from the proposed DSM-PINNs were compared with both existing literature data and alternative computational methods. The BCs considered include simply supported (S-S), cantilever (C-F), fully clamped (C-C), and clamped-simply supported (C-S), where C, F, and S denote clamped, free, and simply supported ends, respectively. The analysis employs the Euler-Bernoulli beam theory, with coordinates illustrated in Fig. \ref{fig4}. For S-S, C-C, and C-S BCs, the material properties are as follows: Young's modulus $E=28$ GPa, density $\rho=2350$ kg/m$^3$, beam length $L=10$ m, width $b=0.2$ m, and height $h=0.6$ m \citep{PatilD1}. For the cantilever beam, the parameters are $E=210$ GPa, $\rho=7860$ kg/m$^3$, $L=0.24$ m, $b=0.012$ m, and $h=0.02$ m. The symbol $\omega_i$ denotes the natural frequency of the $i$-th mode. The DSM-PINNs for eigenvalue problems were optimised employing the Adam optimiser with the $\tanh$ activation function. The deep learning neural network layout comprised two hidden layers, each containing $80$ neurons, denoted as $[1,80,80,2]$. Training was performed on a system including a regular NVIDIA GeForce RTX 5080 GPU.

\begin{figure}[htbp]
  \centering
  \includegraphics[scale=0.8]{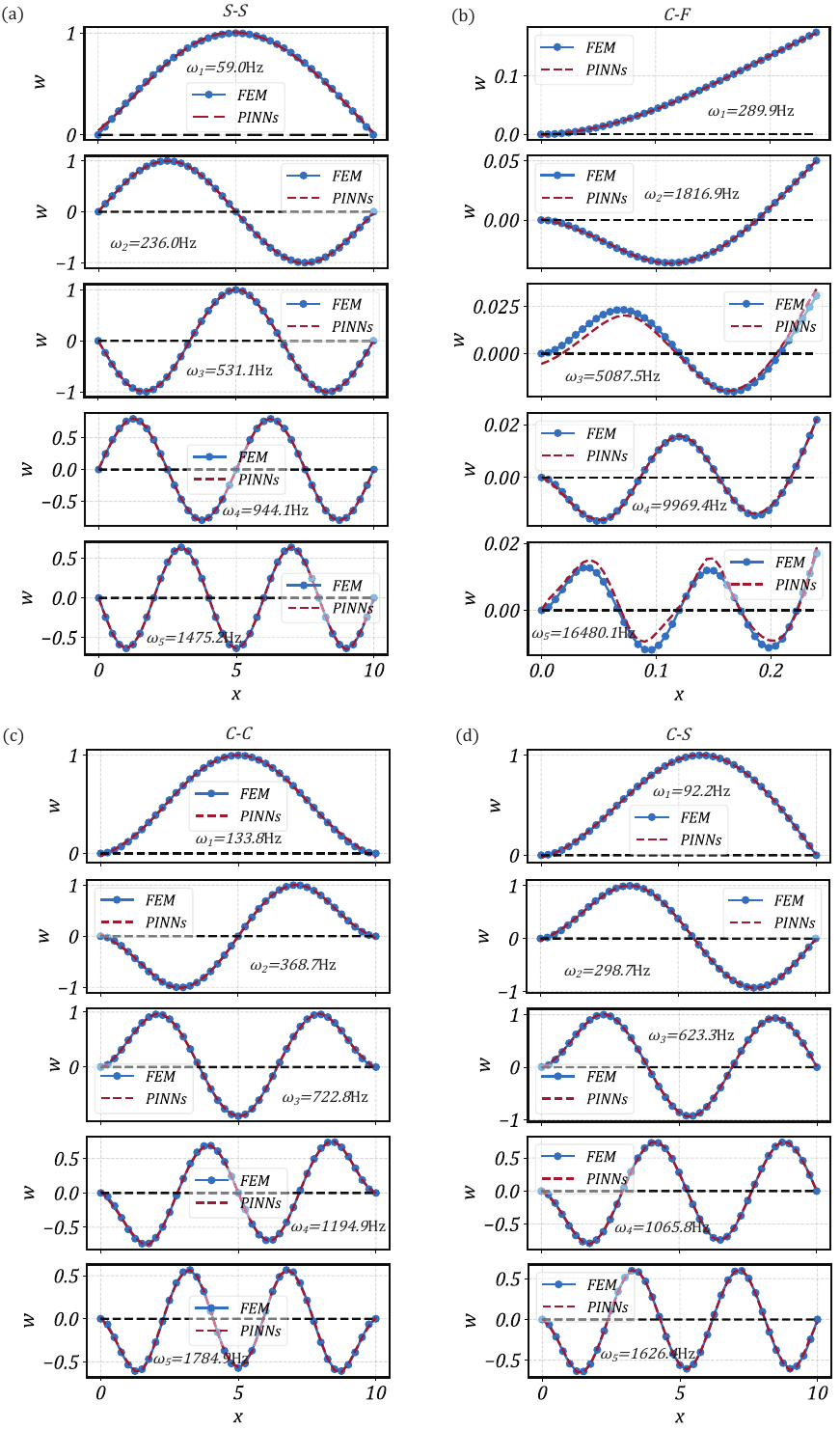}
\caption{The first five eigenfrequencies and the associated mode forms of the beam under distinct boundary conditions (a) S-S, (b) C-F, (c) C-C, and (d) C-S are computed by the proposed PINNs (red dashed line) and FEM (blue line).}
\label{fig5}
\end{figure}

Tabs. \ref{tab1}-\ref{tab2} present the first five natural frequencies of the beam for various BCs. The numerical results are contrasted against results from prior studies utilising the finite element approach \citep{PatilD1}, and DSM using the bisection method. The contrasting results indicate a satisfactory concordance at the identical iteration tolerance ($1e-4$). Then, the first five modeshapes for different BCs are analysed and compared with FEM results, as illustrated in Fig. \ref{fig5}. The modeshape results also agree quite well with those from FEM, except for modes $3$ and $5$ in the cantilever case because the geometric size of the cantilever beam is much lower than that of the other three cases, leading to a relatively higher eigenvalue $\omega_i$. 
\begin{table}[htbp]
  \centering
\scriptsize 
  \caption{The initial five natural frequencies according to Euler-Bernoulli theory with simply supported and cantilever boundary conditions.}
    \begin{tabular}{ccccccc}
    \toprule
   {Mode i} & & Simply-supported (Hz)& &  & Cantilever (Hz)&  \\
    \cmidrule(lr){2-4}
    \cmidrule(lr){5-7}
    & {PINNs} &{DSM} & {\cite{PatilD1}} &{PINNs} & {DSM} & {\cite{PatilD1}}  \\
    \midrule
 1 &  59.0095  & 59.0075  & 58.7520  & 289.9250  & 289.9255  & 289.9253  \\
  2 &  236.0418  & 236.0295  & 227.5096  & 1816.9310  & 1816.9310  & 1816.9311  \\
  3 &  531.0709  & 531.0655  & 488.5903  & 5087.4590  & 5087.4594  & 5087.4593  \\
  4 &  944.1164  & 944.1165  & 818.3652  & 9969.3890  & 9969.3898  & 9969.3898  \\
  5 &  1475.1825  & 1475.1815  & 1196.8516  & 16480.1200  & 16480.1157  & 16480.1158  \\
           \bottomrule
    \end{tabular}%
  \label{tab1}%
\end{table}%

\begin{table}[htbp]
  \centering
\scriptsize 
  \caption{The initial five natural frequencies according to Euler-Bernoulli beam with fully clamped and clamped-simple support boundary conditions.}
    \begin{tabular}{ccccccc}
    \toprule
   {Mode i} & & Fully clamped (Hz)& &  & Clamped-simple support (Hz)&  \\
    \cmidrule(lr){2-4}
    \cmidrule(lr){5-7}
    & {PINNs} &{DSM} & {FEM} &{PINNs} & {DSM} & {FEM}  \\
    \midrule
   1 &133.7631  & 133.7625  & 133.1833  & 92.1899  & 92.1805  & 91.7814  \\
   2& 368.7225  & 368.7225  & 355.4129  & 298.7243  & 298.7245  & 287.9416  \\
   3& 722.8437  & 722.8435  & 665.0297  & 623.2642  & 623.2645  & 573.4152  \\
   4& 1194.8976  & 1194.8965  & 1035.7426  & 1065.8212  & 1065.8185  & 923.8571  \\
   5& 1784.9703  & 1784.9695  & 1448.1904  & 1626.3890  & 1626.3875  & 1319.5288  \\
           \bottomrule
    \end{tabular}%
  \label{tab2}%
\end{table}%
 
\begin{figure}[htbp]
  \centering
  \includegraphics[scale=0.9]{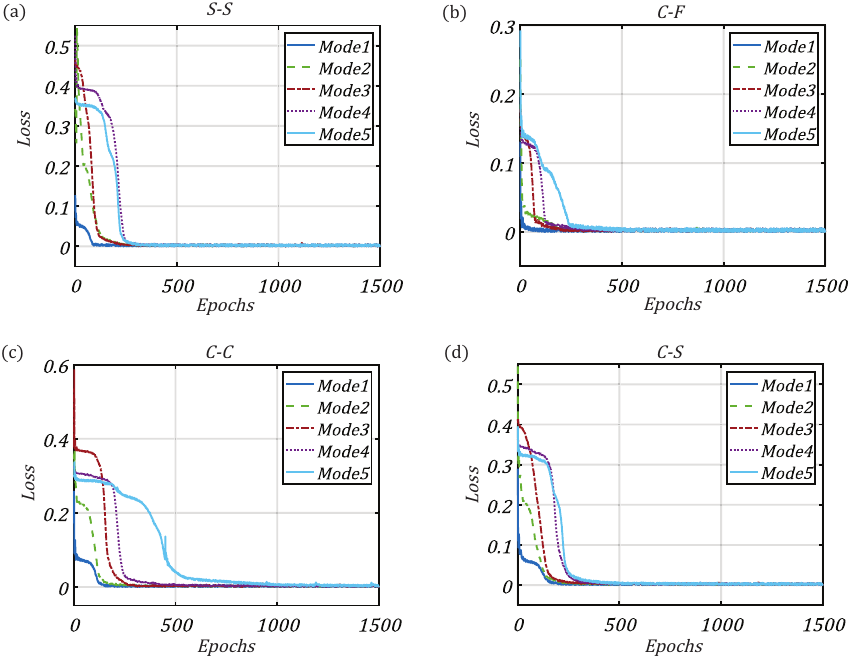}
\caption{The loss of eigenvalue solving for the beam under different boundary conditions while the training procedure is being carried out. (a) simple-supported; (b) cantilever; (c) fully clamped; (d) clamped simple support.}
\label{fig6}
\end{figure}

Fig. \ref{fig6} shows the loss function of the proposed approach during training. The rate of variation in space and time is related to the speed at which deformation waves can travel through the structure. Hence, for higher-order eigenvalues, the neural network takes more iterations to approximate a more tortuous eigenvector, even though the eigenvalue is the `lowest' one. The convergence efficiency is quite high compared with the results reported in \citep{YooS} and the initial drop phase in the loss function for different BCs are completed within 1000 epochs. Moreover, the total number of epochs, including the slow descent phase, is less than 5000.
 
\begin{figure}[htbp]
  \centering
  \includegraphics[scale=0.9]{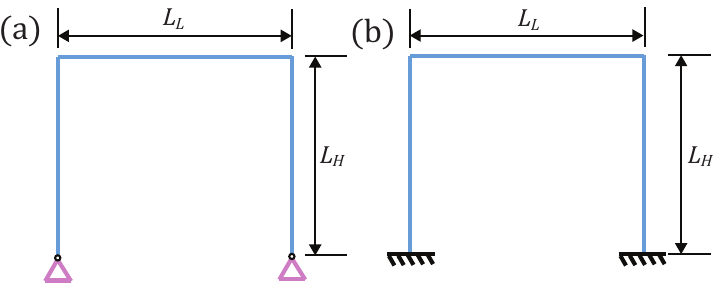}
\caption{Sketch of one bay framed structure. (a) simple-supported; (b) fully clamped.}
\label{fig7}
\end{figure}

Furthermore, two types of frame structures based on Euler-Bernoulli beam theory are analysed using the proposed method to validate the two-dimensional problem. The dimensionless natural frequency $\lambda_i$ is defined to ensure consistent and generalisable outcomes \citep{Banerjee1, ChenZ2}:
\begin{equation}
 \label{eq46} 
\lambda_i=\omega_i\sqrt{\frac{\rho AL^4}{EI}}.
 \end{equation}
The frame structure for a one-bay, one-story system has been studied, as illustrated in Fig. \ref{fig7}, with two boundary conditions: simply supported and clamped. The Young's modulus $E=200$GPa, density $7500$kg/m$^3$ and height $h=0.04$m, and width $b=0.02$m are used. The frame height $L_H$ and beam length $L_L$ are equal to each other $L_H=L_L=1$m. The neural network architecture was $[2,80,80,2]$ because the input variables are 2-dimensional coordinates. The natural frequencies and modeshapes of a frame are obtained. The initial three frequencies are compared with previous works \citep{Banerjee1}, as reported in Tab. \ref{tab3}. 

\begin{table}[htbp]
  \centering
\scriptsize 
  \caption{The initial three dimensionless frequency of the one bay framed structure.}
    \begin{tabular}{ccccccc}
    \toprule
   {Mode i} & & Fully clamped& &  & Simple support&  \\
    \cmidrule(lr){2-4}
    \cmidrule(lr){5-7}
    & {PINNs} &{\cite{Banerjee1}} & {FEM} &{PINNs} & {\cite{Banerjee1}} & {FEM}  \\
    \midrule
   1& 3.2050  & 3.2040  & 3.2050  & 1.4628  & 1.4630  & 1.4630  \\
   2&  12.6391  & 12.6390  & 12.6480  & 9.8664  & 9.8660  & 9.8700  \\
   3& 20.6272  & 20.6270  & 20.6290  & 14.8541  & 14.8540  & 14.8560  \\
           \bottomrule
    \end{tabular}%
  \label{tab3}%
\end{table}%

\begin{figure}[htbp]
  \centering
  \includegraphics[scale=0.85]{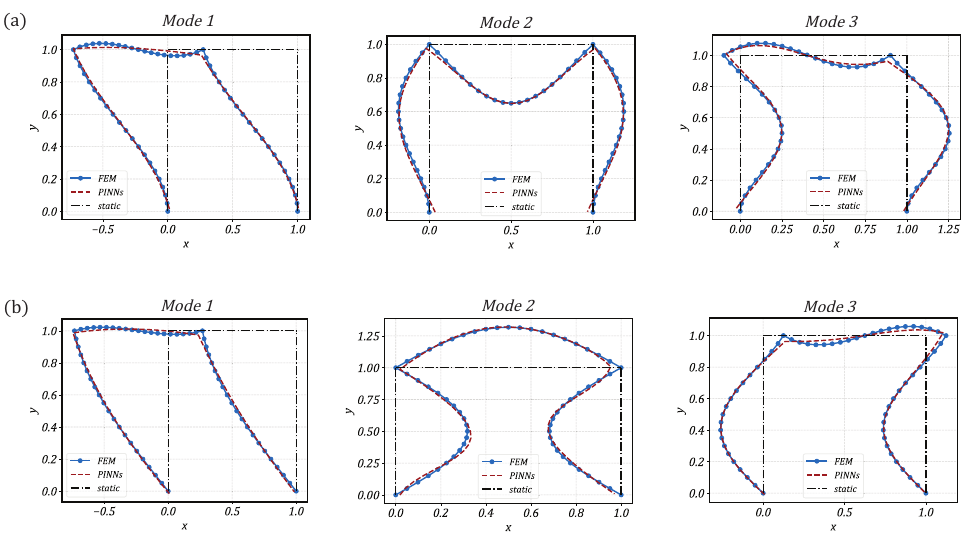}
\caption{First three mode shapes of the framed structures under four different boundary conditions: (a) fully clamped; (b) simple-supported.}
\label{fig8}
\end{figure} 

\begin{figure}[htbp]
  \centering
  \includegraphics[scale=0.9]{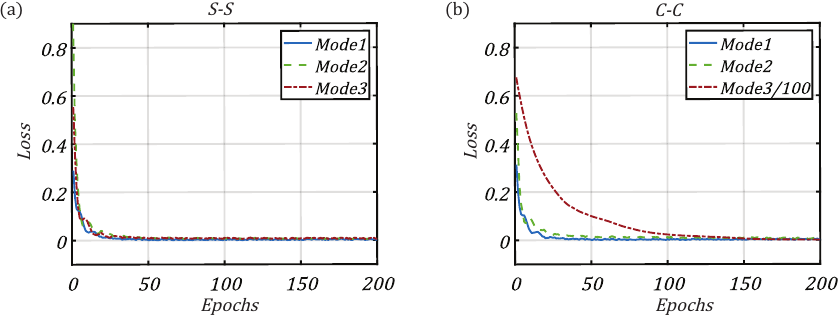}
\caption{The loss of eigenvalue solving for the beam under different boundary conditions throughout the training procedure. (a) simple-supported; (b) fully clamped.}
\label{fig9}
\end{figure}

Table \ref{tab3} demonstrates that the numerical results calculated using the proposed methodology align closely with previous studies and FEM analyses. The modeshapes corresponding to the first three eigenfrequencies, presented in Fig. \ref{fig8}, also exhibit strong agreement with FEM results. Furthermore, the convergence efficiency remains high, as indicated by the loss function evolution in Fig. \ref{fig9}. For simply supported and fully clamped boundary conditions, the iterations converge within 100 and approximately 200 epochs, respectively. Notably, the third mode under fully clamped conditions presents greater computational difficulty compared to the first two modes. Results from both single-beam and frame structures indicate that the proposed DSM-PINNs approach can determine arbitrary-order eigenmodes using the Wittrick-Williams algorithm, achieving acceptable accuracy and convergence efficiency. Additionally, the spectral form of the DSM, specifically, the frequency-domain spectral element method combined with frequency-domain PINNs, which could be employed to determine the dynamic response of the structure.

\section{Addressing the structural dynamics using a frequency domain PINNs}
In the following section, we will investigate the dynamic response of the structure to moving point and impulsive load conditions. The loss function that was presented in Section 2 is chosen in accordance with the characteristics of the dynamic stiffness matrix. Since neural networks are unable to directly process complex numbers, it is necessary to employ two separate networks to train the real and imaginary parts of the problem. Obtaining the dynamic response requires adding the two outputs together and then applying the inverse Fourier transform to the result.

\subsection{Frequency domain PINNs}
\begin{figure}[htbp]
  \centering
  \includegraphics[scale=0.9]{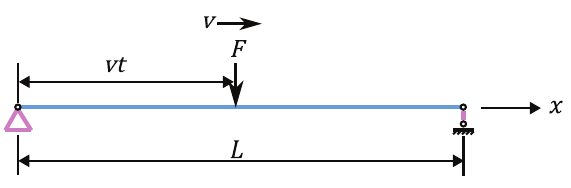}
\caption{Evaluation of a beam subjected to a dynamic point load at a constant velocity.}
\label{fig10}
\end{figure}

Structural dynamics problems frequently involve multiple frequencies, which traditional neural networks may not model with sufficient accuracy. To address this limitation, the present study introduces the frequency-domain (FD) PINN method. The method employs the Fourier transform to translate the original time-domain partial differential equation into the frequency-domain, yielding a DSM and a frequency-domain load. The response spectrum, derived by solving the problem with PINNs in the frequency-domain, is subsequently converted back to the time-domain through the inverse Fourier transform. As an illustrative example, the scenario of a point load moving at constant velocity across an undamped Euler-Bernoulli beam is examined, as depicted in Fig. \ref{fig10}. The corresponding time-domain PDE is:
\begin{equation}
 \label{eq47} 
EI\frac{\partial^4w(x,t)}{\partial x^4}+m\frac{\partial^2w(x,t)}{\partial t^2}=F\delta(x-vt),
 \end{equation} 
where $v$ denotes the velocity of the moving load and $\delta$ indicates the Dirac delta function. If the scenario involves constant acceleration $a$, the right-hand side term becomes $F\delta(x-0.5at^2-v_0t)$.

\begin{figure}[htbp]
  \centering
  \includegraphics[scale=0.9]{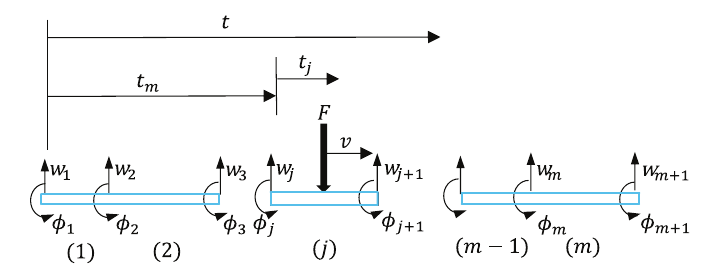}
\caption{Moving load applied to the $j$-th element of the beam.}
\label{fig11}
\end{figure}

The frequency-domain expression of Eq. \eqref{eq47} is derived via the Fourier transform \citep{LiangR1, SarvestanV}:
\begin{equation}
 \label{eq48} 
EI\frac{d^4\bar{w}(x,\omega)}{d x^4}-m\omega^2\bar{w}=\frac{F}{v}e^{\frac{ix\omega}{v}},
 \end{equation} 
where $\bar{w}$ represents the frequency-domain displacements. Recall that this differs from the amplitude obtained via the separation of variables method, but shares the same homogeneous solution introduced in Section 3. Then, the frequency-domain load can be transformed into nodal shear forces and bending moments using the virtual work principle. Using the Fourier transform, the moving point load in the frequency-domain $\textbf{f}$, when the load is in the $j$-th element of the beam as demonstrated in Fig. \ref{fig11}, can be derived: 
\begin{equation}
 \label{eq49} 
\textbf{f}=\frac{F}{v}e^{-i\omega t_m}e^{-i\omega\frac{x_j}{v}},
 \end{equation} 
where $t_m$ is the length of time required for the moving point load to arrive at the $j$-th element (passing over $(j-1)$ elements). $x_j$ is the relative passing distance after time $t_j$, for the constant velocity case, $x_j=t_j*v$.

According to the virtual work principle, a concentrated force applied to an element can be represented as equivalent shear forces and bending moments at the nodes using the shape functions:
\begin{equation}
 \label{eq50} 
\textbf{f}_d=\int^L_0\textbf{N}^T(x,\omega)\frac{F}{v}e^{-i\omega t_m}e^{-i\omega\frac{x_j}{v}}dx.
 \end{equation} 
For the DSM, the shape function is not derived from a polynomial but from a homogeneous solution. Therefore, the equivalent nodal forces $\textbf{f}_d$ are `exact'. The shape function for DSM is:
\begin{equation}
 \label{eq51} 
\bar{w}=[\sin \alpha \zeta, \cos \alpha \zeta, \sinh \alpha \zeta, \cosh \alpha \zeta] \textbf{R}^{-1} \triangle=\textbf{N} \triangle.
 \end{equation} 
Because the matrix $\textbf{R}$ does not include $\zeta=x/L$, the integral \eqref{eq50} can be simplified to:
\begin{equation}
 \label{eq55}
\textbf{f}_d=\frac{F}{v\Delta t}e^{i\omega t_m}(\textbf{R}^{-1})^T\int^1_0[\sin \alpha \zeta , \cos \alpha \zeta, \sinh \alpha \zeta , \cosh \alpha \zeta ]^Te^{i\omega \frac{\zeta}{vL}}d\zeta.
 \end{equation} 
where $\Delta t$ is the time interval for calculation and the moving point load is a function defined on $\Delta t$. However, Eqs. \eqref{eq49} and \eqref{eq50} cannot be directly utilised to the represent moving load and the complete elemental nodal force $\textbf{f}_d$ is shown in the Appendix.

The global nodal force vector $\textbf{P}$ is assembled from $\textbf{f}_d$ using similar methods to those in FEM. Because the stiffness matrix and nodal forces are derived from the governing equation, the system of linear equations can represent the original PDEs:
\begin{equation}
 \label{eq56}
\textbf{W}^*\Delta=\textbf{P}^* \ \Leftrightarrow \ EI\frac{d^4\bar{w}(x,\omega)}{d x^4}-m\omega^2\bar{w}=\frac{F}{v}e^{\frac{ix\omega}{v}}.
 \end{equation}  
 
In addition, because neural networks cannot handle complex numbers, the general displacement can be represented in the form of separated real and imaginary parts, and PINNs solve the above equations:
\begin{equation}
 \label{eq57}
\bar{w}=\bar{w}_r+i\bar{w}_i, \ \ \bar{\phi}=\bar{\phi}_r+i\bar{\phi}_i,
 \end{equation}  
where $\bar{w}_r=\textmd{real}(\bar{w})$ and $\bar{w}_i=\textmd{imag}(\bar{w})$; the same applies to the rotation $\bar{\phi}$. The BCs for the equation can be embedded using methods introduced in Section 2. Therefore, the neural network output $\mathcal{N}(\textbf{x}; \theta)$ can approximate $\bar{w}$ and $\bar{\phi}$ in both real and imaginary parts, thereby formulating our proposed DSM-PINNs as a specific variant of FD-PINNs when solving dynamic response. Accounting for specific initial conditions for the structure in the frequency domain is generally challenging. In most analyses, the response of beams subjected to moving loads assumes a stationary state (null initial displacement and velocity) as the initial condition.

As the network employs frequency inputs, it is essential to delineate the frequency range of the computational domain before computing. Thus, the initial problem is approximated by choosing a frequency range that includes its dominant frequency components. The established criteria for choosing computational frequency ranges in the spectral element method can guide the selection process for DSM-PINNs. To a large extent, the fundamental mode is responsible for determining the dynamic response of simply supported beams when subjected to moving loads \citep{LiangR1}, so the frequency input range used by the neural network's input needs to include the natural frequency of the first order of the structure. In all other circumstances, the computational frequency range needs to be determined by taking into account the frequency components of the load being applied. The eigenmode can be derived using the techniques outlined in the preceding section.

\begin{figure}[htbp]
  \centering
  \includegraphics[scale=0.8]{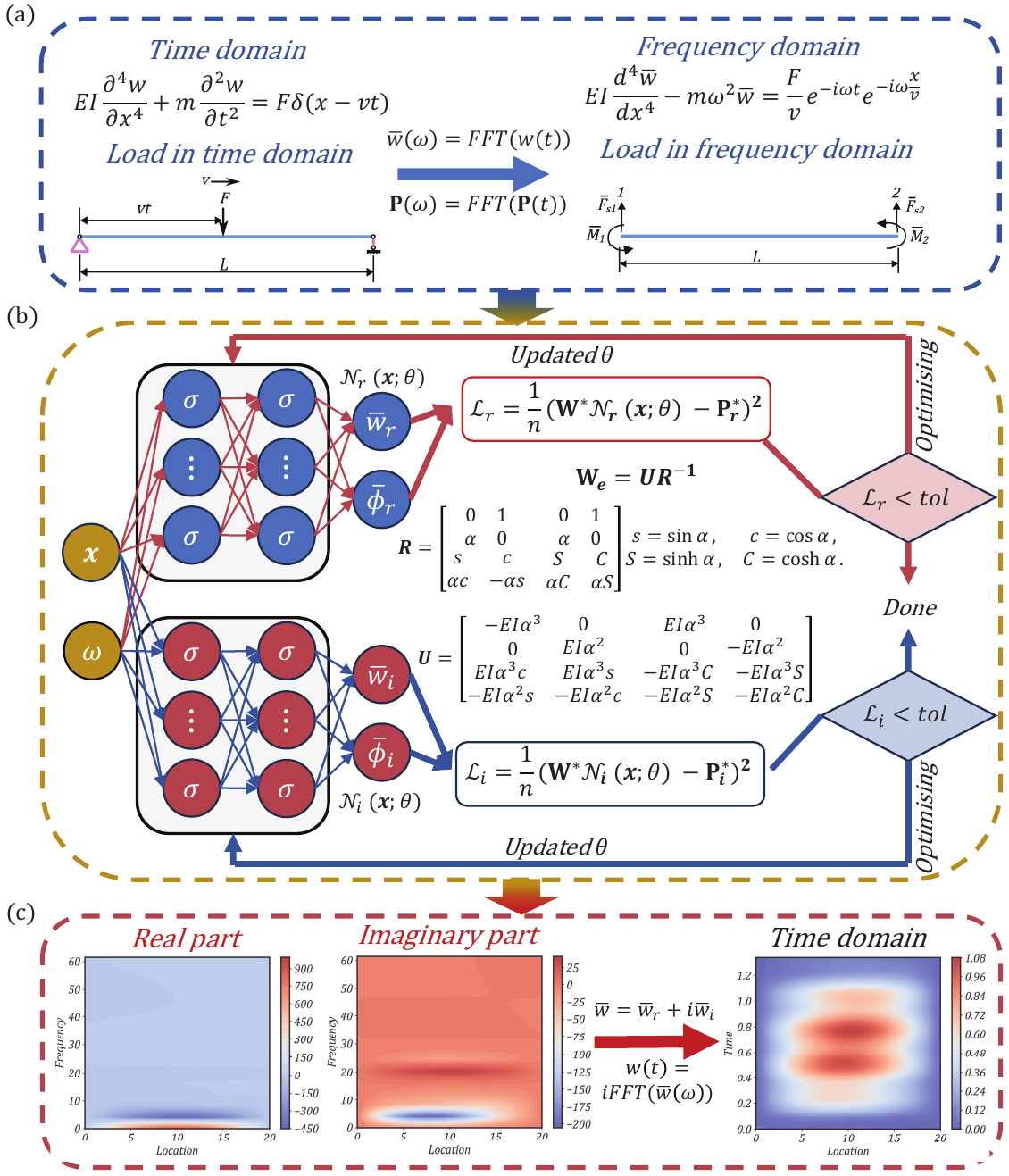}
\caption{Methodology for addressing structural dynamic response with the DSM-PINNs approach. (a). The governing equation and dynamic loading are transformation into the frequency domain, exemplified by a moving point load traversing the Euler-Bernoulli beam. (b) In determining the frequency-domain dynamic response via FD-PINNs, the nodal equilibrium equation is selected as the loss function. Two distinct neural networks are employed for the real and imaginary aspects. (c). Through the application of the inverse Fourier transform, the time-domain dynamic response can be acquired. This method can be used to other structural type elements and dynamic loads.}
\label{fig12}
\end{figure}

Following the completion of the FD-PINN calculation for the beam's response spectrum in the frequency domain, the discrete inverse Fourier transform is applied to determine the solution of the governing equations in the time domain, as shown below:
\begin{equation}
 \label{eq58}
w(x,t_n)=iFFT(\bar{w}_r(x,\omega_k)+i\bar{w}_i(x,\omega_k))=\frac{1}{N}\sum^{N-1}_{k=0}\bar{w}(x,\omega_k)e^{\frac{2ikn\pi}{N}},
 \end{equation}  
where $\bar{w}(x,\omega_k)$ is the $k$-th component in the frequency domain and composed of two parts $\bar{w}_r(x,\omega_k)$ and $i\bar{w}_i(x,\omega_k)$, $w(x,t_n)$ is the $n$-th associated displacement component in the time domain. Furthermore, the data from measurements taken in the time domain can be converted into the frequency domain and incorporated into the loss function as outlined in Eq. \eqref{eq14} via the Fourier transform. The complete procedure for solving a beam under moving point loads using the DSM-PINNs is shown in Fig. \ref{fig12}.

\subsection{The properties of dynamic stiffness matrix and loss function selected}
The $\textbf{K}$ stiffness matrix in traditional FEM is a positive semi definite matrix without restrains and a positive definite matrix $\textbf{K}^*$ after applying sufficient BCs, which means for an arbitrary vector $\textbf{y}$:
\begin{equation}
 \label{eq59}
\textbf{y}^T\textbf{K}\textbf{y}\geq0, \ \ \textbf{y}^T\textbf{K}^*\textbf{y}>0.
 \end{equation}  

Consequently, the first-order derivative of the loss function is zero, and the Hessian matrix is positive definite in both Eqs. \eqref{eq20} and \eqref{eq23}, resulting in the neural networks converging to a local minimum. This indicates that these methods exhibit a certain degree of local convergence stability. In contrast, the DSM $\textbf{W}(\omega)$ contains elements that are functions of the trial frequency $\omega$, rendering the matrix indefinite depending on the value of $\omega$. The following is an example of how the DSM can be expanded using a first-order Taylor series:
\begin{equation}
 \label{eq60}
\textbf{W}(\omega)\approx \textbf{K}-\omega^2\textbf{M}+o(\omega^2),
 \end{equation}  
where $K$ and $M$ are the stiffness matrix and the associated mass matrix (also called the consistent mass matrix) in standard FEM, the properties of the dynamic stiffness matrix $\textbf{W}(\omega)$ can be represented by the matrix $ \textbf{K}- \omega^2 \textbf{M}$. Assuming the $\lambda$ and $\textbf{q}$ are the standard eigenvalue and eigenvector of matrix $\textbf{K}-\omega^2\textbf{M}$:
\begin{equation}
 \label{eq61}
(\textbf{K}-\omega^2\textbf{M})\textbf{q}=\lambda \textbf{q}.
 \end{equation}  
The properties of the matrix $\textbf{K}-\omega^2\textbf{M}$ can be encapsulated in the following form:
\par
(1) if the trial frequency $\omega^*$ is the general eigenvalue of $\textbf{K}-{\omega^*}^2\textbf{M}$, indicating that the determinant $|\textbf{K}-{\omega^*}^2\textbf{M}|=0$. The case where the standard eigenvalue $\lambda^*=0$ yields the modeshpae eigenvector $\textbf{q}^*$ corresponding to $\omega^*$;
\par
(2) if the trial frequency $\omega^*$ is lower than the lowest general eigenvalue (first order mode) $\omega_1$, so that $\omega^*<\omega_1$. The dynamic stiffness matrix $(\textbf{K}-{\omega^*}^2\textbf{M})$ is a positive definite matrix after sufficient BCs;
\par
(3) if the trial frequency $\omega^*$ is higher than $i$-th order general eigenvalue $\omega_i$, so that $\omega^*>\omega_i$. The stiffness matrix $(\textbf{K}-{\omega^*}^2\textbf{M})$ is an indefinite matrix, and has $i$ number negative standard eigenvalue $\lambda$;
\par
These properties can be proved using the generalized eigenvalues and eigenvectors via their orthogonality with respect to the matrix $(\textbf{K}-\omega^2\textbf{M})$. Note that these properties can also be used in FEM dynamics problems.

Therefore, the energy-type loss function cannot be used if the DSM is employed. For the nodal equilibrium equations-based PINNs, $(\textbf{K}-\omega^2\textbf{M})^2$ is a semi-positive definite matrix ($|\textbf{W}(\omega)|=0$ as $\omega$ associated with a natural frequency), so the resulting Hessian matrix is also positive semi-definite \eqref{eq23}. If the appropriate trial frequency is chosen, such that the second-order derivative-Hessian matrix is also positive definite. This type of loss function can be used in DSM-PINNs, leading to a strong-form representation of the original PDEs without energy variation:
\begin{equation}
 \label{eq62}
\mathcal{L}_{phys}=\frac{1}{n}\sum^n_{i=1}\left[\textbf{K}_{ij}\cdot{\mathcal{N}(\textbf{x}_j; \theta)}-\textbf{P}_i\right]^2.
 \end{equation} 
 
\subsection{The dynamic response for moving point load passing structure}
The performance of DSM-PINNs in solving structural dynamic problems is evaluated through three case studies. Results obtained from the software Abaqus using B23 element, consistent mass matrix and total $40$ elements in a single beam with implicit time-integration serve as a reference with time steps of $0.016$, $0.0125$ and $0.005$ for the three cases, as replicated in previous work \citep{LiangR1}. The beam length is $20$m, density per unit length $3000$kg/m, and stiffness $EI=2.4\times10^9$N$\cdot$m$^2$. The velocities of the moving point load are $15$m/s, $20$m/s, and $50$m/s, each with a load value of $12$kN. For each case, a total of $N=256$ are used at a selected point on the structure, and the time steps are $0.006$s, $0.004$s, and $0.002$s, respectively. Based on the discrete Fourier transform, the frequency resolutions are $4.091$rad/s, $6.1359$rad/s, and $12.272$rad/s ($\Delta\omega=2\pi/(N\Delta t)$), ensuring that the full frequency range encompasses the fundamental natural frequency of the beam. The elements in the dynamic stiffness matrix are entirely real, while only the force is expressed as a complex number. The real and imaginary components (displacements) are trained independently. All results are presented in dimensionless form using the scaling factor $FL^3/(48EI)$.

\begin{figure}[htbp]
  \centering
  \includegraphics[scale=0.9]{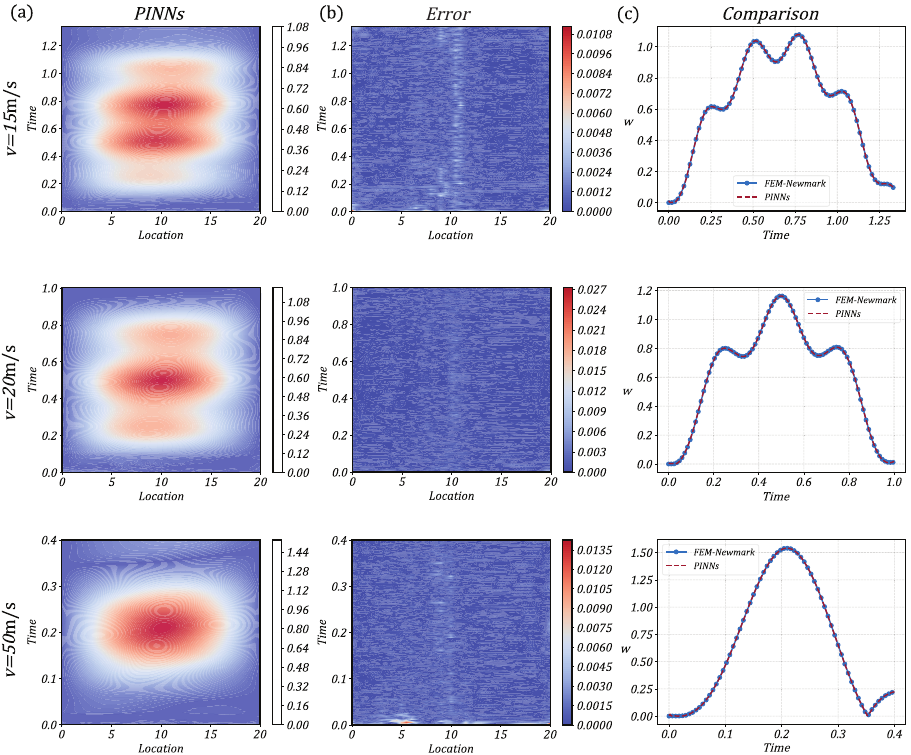}
\caption{The dimensionless displacement field was calculated by PINNs with three different moving load velocities. (a). The results from DSM-PINNs. (b). The error compared with the results from FEM. (c). The results comparison at the fixed midspan location.}
\label{fig13}
\end{figure}

The DSM-PINNs developed for dynamic response problems were trained employing the Adam optimiser, with the $tanh$ activation function. The architecture of the deep learaning model (neural network) was $[2, 200, 200, 200, 2]$, comprising three hidden layers, each containing 200 neurons. Training was conducted on a standard NVIDIA GeForce RTX 5080 GPU. Fig. \ref{fig13} displays the displacement outcomes computed by the fully optimised model across different situations. For systems without damping, it is unnecessary to distinguish between transient and steady-state responses. The relative errors are calculated as follows:
\begin{equation}
 \label{eq63}
Error=\frac{(w_{PINNs}-w_{FEM})}{w_{FEM}}, 
 \end{equation} 
where $w_{PINNs}$ is the displacement obtained from DSM-PINNs, and $w_{FEM}$ is the displacement obtained from FEM. 

\begin{figure}[htbp]
  \centering
  \includegraphics[scale=0.9]{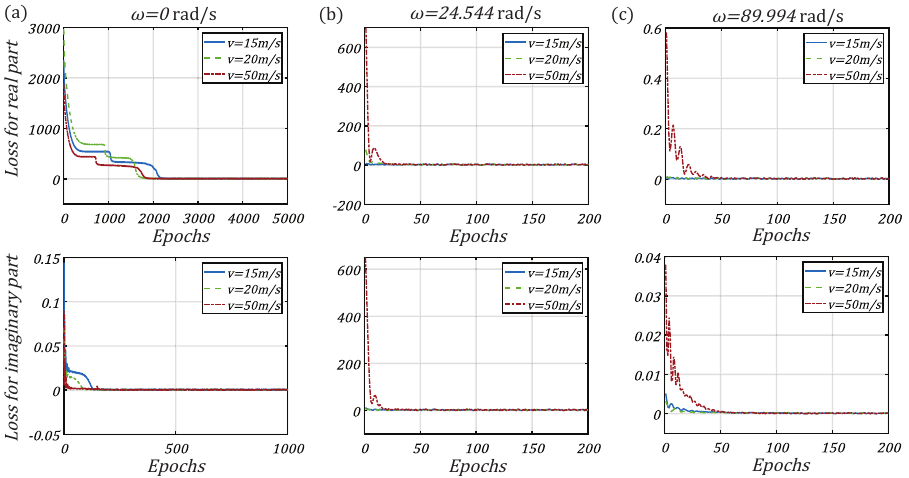}
\caption{The value for the loss function during model training at different trial frequencies for both real and imaginary parts: (a). at $\omega=0$rad/s. (b) at $\omega=24.544$rad/s. (c) at $\omega=89.9644$rad/s.}
\label{figmovingloss}
\end{figure}

The results obtained from the proposed DSM-PINNs demonstrate low and acceptable relative errors compared to the finite element method (FEM) across all three cases, with a maximum error of approximately $0.027$. Based on this, it can be concluded that the method is capable of accurately simulating the dynamic response of beams when subjected to moving loads. This variation in the loss function of the developed models is shown in Fig. \ref{figmovingloss}, which illustrates the relationship between the number of training iterations and the loss function. The calculated natural frequencies are $\omega_1=22.112$ rad/s and $\omega_2=88.275$ rad/s, and the selected trial frequencies are $\omega=0$ rad/s, $\omega=24.544$ rad/s, and $\omega=89.964$ rad/s. The model exhibits efficient convergence, typically within $100$ iterations, except for the case with $\omega=0$ rad/s for the real part, which corresponds to a static scenario, as evidenced by the negligible loss for the imaginary part. In addition, increasing the velocity of the moving load has a tendency to produce an increase in the number of iterations that are necessary for the model to converge.

Therefore, the neural network is no longer obligated to replicate the original complex function with multi-frequency characteristics, thereby significantly enhancing the simulation accuracy.

\subsection{The dynamic response for a unit impulsive point force on structure}
\begin{figure}[htbp]
  \centering
  \includegraphics[scale=0.7]{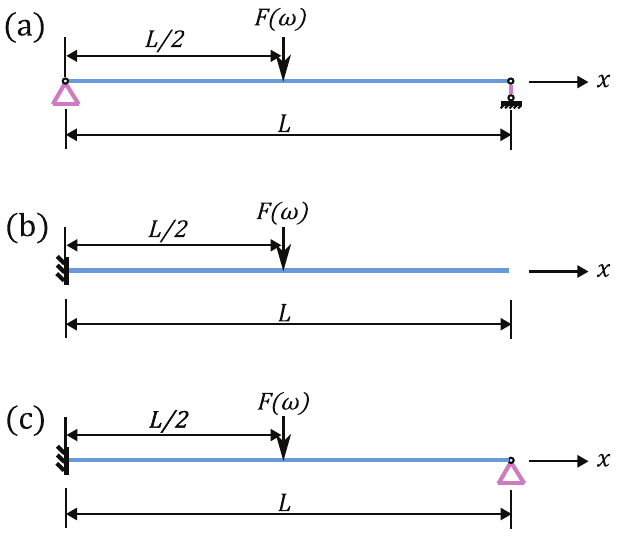}
\caption{An Euler-Bernoulli beam that is exposed to a transverse unit impulsive point load while being subjected to a variety of boundary conditions: (a) simply supported, (b) cantilever, (c) clamped-simply supported.}
\label{fig14}
\end{figure}

An examination of the dynamic response of a beam that has been exposed to an impulsive point load under a variety of boundary conditions is presented in this subsection, as illustrated in Fig. \ref{fig14}, consistent with previous studies employing the spectral element method \citep{KimT2}. Results obtained from the software Abaqus using B23 element, consistent mass matrix and total $40$ elements in a single beam with implicit time-integration serve as a reference. Both the geometric and material properties of the example beam are identical to those that were used in the case of the moving point load simulation, specifically $L=20$ m, $m=3000$ kg/m, and $EI=2.4\times10^9$ N$\cdot$m$^2$. The architecture of the deep learning model employed was $[2, 200, 200, 200, 2]$ with a $tanh$ activation function, and training was conducted using the Adam optimiser.

\begin{figure}[htbp]
  \centering
  \includegraphics[scale=0.9]{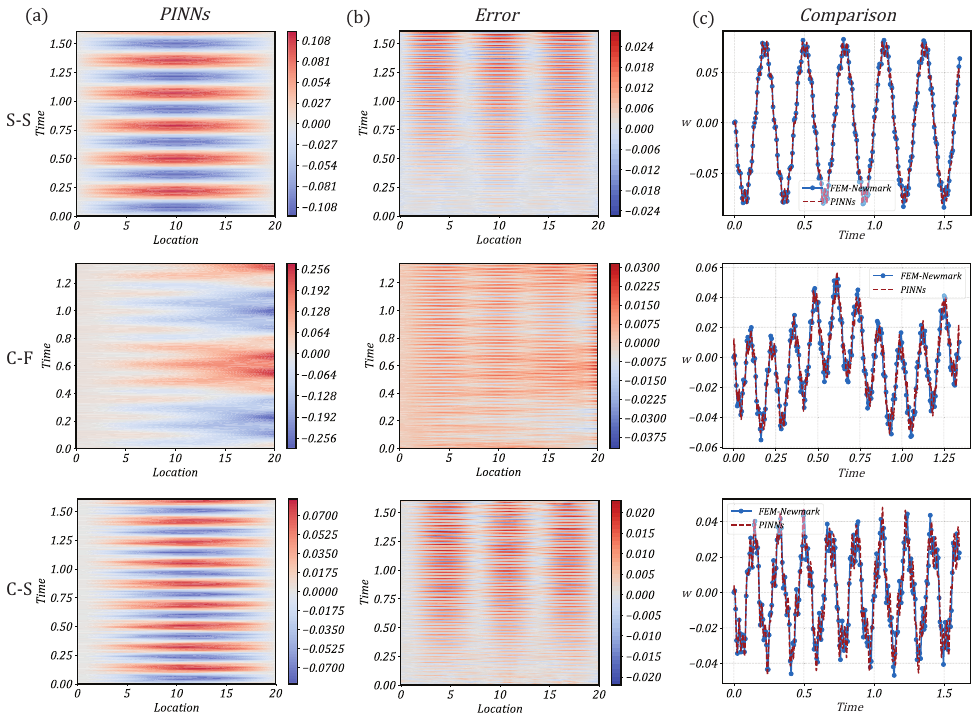}
\caption{An Euler-Bernoulli beam that is exposed to a transverse unit impulsive point load while being subjected to a variety of boundary conditions: (a) the results from DSM-PINNs, (b) the error compared with results from FEM, (c) the results comparison at $x/L=0.25$ location.}
\label{fig15}
\end{figure}

\begin{figure}[htbp]
  \centering
  \includegraphics[scale=0.9]{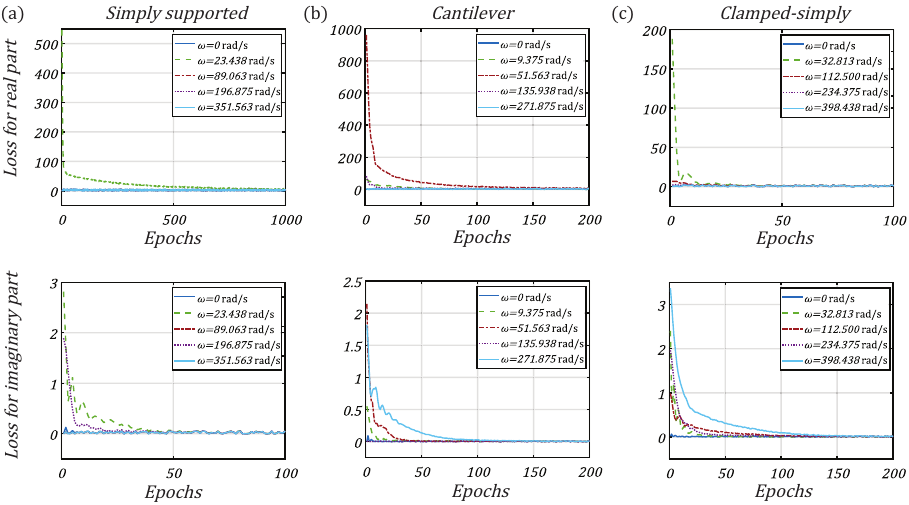}
\caption{The value for the loss function during model training at different trial frequencies for both real and imaginary parts under different boundary conditions: (a). Simply-supported. (b) Cantilever. (c) Clamped-simply supported.}
\label{figlossimpul}
\end{figure}

The dynamic responses are compared with those from the implicit time integration in FEM, and the results for the transverse displacement from the trained model are shown in Fig. \ref{fig15}. Nevertheless, the DSM-PINNs model shows relatively low error (in three cases, the maximum relative error is almost $0.0375$) compared to FEM results. Consequently, the recommended approach results in a neural network that is freed from the responsibility of simulating the original complex function, which possesses multi-frequency characteristics. Thus, DSM-PINNs enhances simulation accuracy and maintains low errors. The evolution of the loss function during model training is illustrated in Fig. \ref{figlossimpul}.

Similar to the moving-load-induced beam problem, the observed frequency point may be near the natural frequency or at $\omega=0$. A calculation is made using the method described in the previous section to determine the first four natural frequencies of the beam under various boundary conditions. The results of this calculation are presented in Tab. \ref{tab4}.
\begin{table}[htbp]
  \centering
\scriptsize 
  \caption{The initial four circular natural frequency of beam structure under different boundary condition.}
    \begin{tabular}{cccc}
    \toprule
   {Mode i} & Simple-supported  & Cantilever& Clamped-simply \\
    \midrule
   1& 22.111  & 7.884  & 34.490  \\
   2&  88.277  & 49.231  & 111.737  \\
   3& 198.623  & 137.960  & 233.107  \\
      4& 353.104  & 270.402  & 398.623  \\
   \bottomrule
    \end{tabular}%
  \label{tab4}%
\end{table}%

As shown in Fig. \ref{figlossimpul}, the model exhibits relatively high convergence efficiency, typically requiring fewer than $200$ iterations, except for cases with $\omega=0$ rad/s for the real part. Moreover, for the real part in all three boundary conditions, the trial frequency closest to the first-order natural frequency takes more epochs to converge. However, the imaginary part corresponding to the fourth-order natural frequency requires more iterations, a phenomenon that warrants further investigation.

\section{Conclusions}
This paper presents DSM-PINNs, a novel method for solving original partial differential equations, whether homogeneous (vibration eigenvalue problems) or nonhomogeneous (dynamic response problems). The context of this paper is restricted to the Euler-Bernoulli beam problem, but can be extended to other structures because the core step is only associated with the dynamic stiffness matrix $\textbf{W}$, not limited to beams. This method offers a framework that is both flexible and robust, integrating physical laws with data gathered from observations. Drawing inspiration from FEM and physics-informed neural networks, the proposed dynamic stiffness matrix is constructed by using homogeneous solutions of partial differential equations rather than polynomial formulations in standard FEM. This method, when combined with PINNs, makes use of an exact shape function for each element, making it possible to formulate the stiffness matrix without energy-based approaches, ultimately resulting in a robust strong-form loss function. A brief summary of the most important findings is as follows:

(1) In the context of the structural eigenvalue problem, there are two primary challenges that are encountered: the first is the acquisition of nontrivial solutions, and the second is the resolution of the tendency of neural networks to converge to the lowest-order eigenmode. By utilising the Wittrick-Williams algorithm, which successfully restricts the eigenvalue within narrow bounds and transforms arbitrary-order eigenvalues into the lowest-order form, these issues are effectively addressed. It is therefore possible for PINNs to converge and generate the mode shapes that correspond appropriately. Given that there is an infinite number of eigenfunction combinations that can satisfy a particular PDE, normalization techniques are utilised during each neural network training session in order to guarantee the existence of nontrivial approaches to solving the problem. Strong agreement can be seen between the results obtained from FEM and those obtained for single-beam and frame structures. Additionally, the iteration checks of the loss function demonstrate that the convergence remains highly efficient. It takes approximately $1000$ epochs for the tracing process to be completed for an individual beam with a variety of boundary conditions, and it takes approximately $200$ epochs for a frame-type structure.

(2) It is possible to think of the DSM as being analogous to the frequency-domain spectral element when discussing structural dynamic response evaluation. Therefore, FD-PINNs are able to solve the problem of spectral bias that arises in neural networks when they are used to simulate multi-frequency functions. The modeling complexity of the neural network is reduced when the displacement field and load are transformed into the frequency domain. Differential terms are incorporated into the DSM, which eliminates the requirement for automatic differentiation. Additionally, the system of linear equations represented at each frequency is a complete representation of the PDEs that were initially used. The results of the calculations for a moving point load traversing a structure and an impulsive load are validated by employing the Newmark implicit time integration method within FEM. It is clear that the proposed DSM-PINNs is both effective and efficient, as evidenced by the fact that the relative error between these results is small. This approach demonstrates significant potential for applications in the fields of computational mechanics, digital twin technologies, and artificial intelligence.

The method can be applied to structures that are more complicated, such as plates and shells, if the DSM $\textbf{W}$ is determined. Particularly noteworthy is the fact that the proposed method is compatible with standard FEM-PINNs. This is because the DSM and FEM share similarities in terms of how they deal with eigenvalue problems and dynamic responses (in fact, the Wittrick-Williams algorithm could also be used in FEM problems). Additionally, the method can be incorporated into deep operator neural networks in order to facilitate the treatment of random boundary conditions and the assembly of such structures. This is accomplished by incorporating families of DSMs within neural networks. In addition, it is possible to develop a unified model that addresses both forward and inverse problems.

\textbf{Declaration of competing interest}:\par
The authors declare that they have no known competing financial interests or personal relationships that could have appeared to influence the work reported in this paper.

\textbf{Acknowledgments}:\par
This work was supported by the National Natural Science Foundation of China (12402135) and the Fundamental Research Funds for the Central Universities.
%

\section*{Appendix}
According to the homogeneous solution of Euler-Bernoulli beam theory \eqref{eq27}, the matrix equations \eqref{eq29} and \eqref{eq30} can be explicitly represented by:
 \begin{equation*}
\Delta=\textbf{R}\textbf{A}=\left[\begin{matrix}\bar{w}_1 \\\bar{\phi}_1 \\ \bar{w}_2 \\ \bar{\phi}_2\\ \end{matrix}\right]=\left[\begin{matrix}0 & 1& 0 & 1 \\\alpha &0 &\alpha &0  \\ s&c & S & C\\ \alpha c & -\alpha s &  \alpha C&\alpha S \\\end{matrix}\right] \ \left[\begin{matrix}A_1 \\A_2 \\ A_3 \\ A_4\\\end{matrix}\right],
 \end{equation*} 
 \begin{equation*}
 \textbf{P}=\textbf{U}\textbf{A}=\left[\begin{matrix}-\bar{Q}_1 \\\bar{M}_1 \\ \bar{Q}_2 \\ -\bar{M}_2\\ \end{matrix}\right]=\left[\begin{matrix}-EI\alpha^3 & 0& EI\alpha^3& 0 \\0 &EI\alpha^2 &0 & -EI\alpha^2  \\ EI\alpha^3c& EI\alpha^3s& -EI\alpha^3C & -EI\alpha^3S \\ -EI\alpha^2s & -EI\alpha^2c& -EI\alpha^2S & -EI\alpha^2C \\\end{matrix}\right] \ \left[\begin{matrix}A_1 \\A_2 \\ A_3 \\ A_4\\\end{matrix}\right],
 \end{equation*}
where $c=\cos \alpha$, $s=\sin \alpha$, $C=\cosh \alpha$ and $S=\sinh \alpha$. Hence, the DSM for the Euler-Bernoulli beam is:
 \begin{equation*}
\textbf{W}(\omega)=\frac{EI}{L^2}\left[\begin{matrix}\delta/L & \theta& -\delta_l/L& \theta_l \\\theta &\beta L &-\theta_l &\gamma L \\ -\delta_l/L& -\theta_l& \delta/L& -\theta \\ \theta_l & \gamma L& -\theta & \beta L \\\end{matrix}\right],
 \end{equation*}
where parameters $\delta=(sC+cS)\alpha^3/D$, $\delta_l=(s+S)\alpha^3/D$, $\theta=sS\alpha^2/D$, $\theta_l=(C-c)\alpha^2/D$, $\beta=(sC-cS)\alpha/D$, $\gamma=(S-s)\alpha/D$ and $D=1-cC$. It is possible to expand the stiffness matrix using a Taylor series to first order \eqref{eq60}:
 \begin{equation*}
\textbf{W}(\omega)\approx \textbf{K}-\omega^2\textbf{M},
 \end{equation*}
 where the matrices $\textbf{K}$ and $\textbf{M}$ are the stiffness and mass matrices in FEM:
 \begin{equation*}
\textbf{K}=\frac{EI}{L^3}\left[\begin{matrix}12 & 6L& -12& 6L \\6L &4L^2 &-6L & 2L^2\\ -12& -6L& 12& -6L \\ 6L & 2L^2& -6L & 4L^2 \\\end{matrix}\right], \ \ \textbf{M}=\frac{mL}{420}\left[\begin{matrix}156 & 22L& 54& -13L \\22L &4L^2 &13L &-3L^2\\ 54& 13L& 156& -22L \\ -13L & -3L^2& -22L & 4 L^2 \\\end{matrix}\right].
 \end{equation*}
The variable $j_0$ denotes the number of natural frequencies of the structure in the range from $\omega=0$ to $\omega=\omega^*$, under the condition that all displacements are restrained.
\begin{equation*}
j_0=\Sigma j_m,
 \end{equation*} 
where $j_m$ denotes the count of natural frequencies associated with each single element of the entire structure, under the assumption that all boundary conditions are fixed:
\begin{equation*} 
j_m=\textmd{highest}\ \textmd{integer}[\frac{\alpha}{\pi}]+0.5\left(2-\textmd{sign}(\beta L)-\textmd{sign}(\beta L -(\gamma L)^2/(\beta L))\right),
 \end{equation*} 

For a two-dimensional frame structure, the frame stiffness matrix is augmented by the axial dynamic stiffness part:
 \begin{equation*}
\textbf{W}_A(\omega)=\frac{EA}{L}\left[\begin{matrix}\epsilon \cot \epsilon & -\epsilon \csc \epsilon\\  -\epsilon \csc \epsilon& \epsilon \cot \epsilon \\\end{matrix}\right],
 \end{equation*}
where $\epsilon=\sqrt{m\omega^2L^2/(EA)}$, which is associated with the longitudinal wave phase speed.

The nodal force $\textbf{f}_d$ as a moving point load passing along the Euler-Bernoulli beam \eqref{eq55} can be further simplified:
 \begin{equation*}
\textbf{f}_d=\frac{F}{v\Delta t}e^{i\omega t_m}(\textbf{R}^{-1})^Tf_d,
 \end{equation*}
where $f_d$ is a $4\times1$ vector given by:
 \begin{equation*}
f_d(1)=\frac{e^{\frac{-iL\omega}{v}}v(e^{\frac{iL\omega}{v}}v\alpha-v\alpha\cos\alpha-iL\omega\sin\alpha)}{v^2\alpha^2-L^2\omega^2}, 
 \end{equation*}
 \begin{equation*} 
 f_d(2)=\frac{e^{\frac{-iL\omega}{v}}v(ie^{\frac{iL\omega}{v}}L\omega-iL\omega\cos\alpha+v\alpha\sin\alpha)}{v^2\alpha^2-L^2\omega^2},
 \end{equation*} 
 \begin{equation*}
f_d(3)=\frac{e^{\frac{-iL\omega}{v}}v(e^{\frac{iL\omega}{v}}v\alpha+v\alpha\cosh\alpha+iL\omega\sinh\alpha)}{v^2\alpha^2+L^2\omega^2},
 \end{equation*}
 \begin{equation*}
 f_d(4)=\frac{e^{\frac{-iL\omega}{v}}v(-ie^{\frac{iL\omega}{v}}L\omega+iL\omega\cosh\alpha+v\alpha\sinh\alpha)}{v^2\alpha^2+L^2\omega^2}.
 \end{equation*}
 
\bibliography{referencesnew} 
\bibliographystyle{elsarticle-harv}

\end{document}